\documentclass[final,11pt]{siamltex}

\usepackage{amsmath,amssymb,amsxtra,comment,graphicx,psfrag}
\usepackage{xspace}
\usepackage{pdfsync}
\usepackage{color}
\usepackage{hyperref}
\usepackage{graphicx}
\usepackage{amssymb}
\usepackage{epstopdf}
\usepackage{subfigure}
\newtheorem{remark}[theorem]{Remark}
\newcommand{\EOC}{\ensuremath{\operatorname{EOC}}}
\newcommand{\IEI}{\ensuremath{\operatorname{IEI}}}

\newcommand{\norm}[2]{\ensuremath{\| #1 \|_{#2}}}

\newcommand{\mean}[1]{\{ \kern -1.6mm \{#1\} \kern -1.6mm \}}
\newcommand{\ha}{\frac{1}{2}}
\newcommand{\jump}[1]{[ \kern -.7mm [#1] \kern -.7mm ]}

\newcommand{\ud}{\mathrm{d}}
\newcommand{\ndg}[1]{| \kern -.25mm \|{#1}| \kern -.25mm \|}
\newcommand{\ltwo}[2]{\|{#1}\|_{#2}}
\newcommand{\ltwoin}[2]{\langle{#1},{#2}\rangle}
\newcommand{\enin}[2]{\langle\!\langle{#1},{#2}\rangle\!\rangle}
\newcommand{\enorm}[1]{\|{#1}\|_a}

\newcommand{\mA}{\mathcal{A}}
\newcommand{\mR}{\mathcal{R}}

\def\dom{\Omega}
\def\timestep{k }
\def\polydegree{p}

\def\initpdesol{u_0}
\def\initpdevel{v_0}

\def\pdesol{u}

\def\usoln{U}
\def\urecsoln{\hat{U}}

\def\velocsoln{V}
\def\recvelocsoln{\hat{V}}

\def\pdesource{f}
\def\pdesourcetimeapprox{\bar{f}}

\renewcommand{\norm}[2]{||#1||_{#2}}

\renewcommand{\enorm}[1]{|||#1|||}

\newcommand{\be}{\begin{equation}}
\newcommand{\ee}[1]{\label{#1} \end{equation}}
\newcommand{\bea}{\begin{eqnarray}}
\newcommand{\eea}{\end{eqnarray}}
\newcommand{\bale}[2]{\begin{equation} \left#1 \begin{array}{#2}}
\newcommand{\eale}[2]{\end{array} \right#1 \label{#2} \end{equation} }

\newcommand{\bal}[2]{\left#1 \begin{array}{#2}}
\newcommand{\eal}[1]{\end{array} \right#1  }

\newcommand{\sobspace}[1]{\operatorname H^{#1}}

\newcommand{\sob}[2]{\sobspace{#1}(#2)}

\newcommand{\sobu}[3]{\sobspace{#1}_{#2}(#3)}

\newcommand{\lebspace}[1]{\operatorname L^{#1}}
\newcommand{\leb}[2]{\lebspace{#1}(#2)}

\newcommand{\email}[1]{Email: \href{mailto:#1}{\texttt{#1}}}
\title{A posteriori error estimates \\ for leap-frog  and cosine methods \\ for second order evolution problems}
\author{
  Emmanuil H.~Georgoulis
  \thanks{
    Department of Mathematics, University of Leicester, Leicester LE1 7RH, England UK, and School of Applied Mathematical and Physical Sciences, National Technical University of Athens, Athens 15780, Greece. 
    \email{Emmanuil.Georgoulis@le.ac.uk}
  }
  \and
  Omar Lakkis
  \thanks{
    Department of Mathematics, University of Sussex, Brighton BN1 9QH, England UK.
    \url{http://www.maths.sussex.ac.uk/Staff/OL}
  }
  \and
  Charalambos G.~Makridakis
  \thanks{
    Department of Mathematics, University of Sussex, Brighton BN1 9QH, England UK.
    \email{C.Makridakis@sussex.ac.uk}
  }
  \and
  Juha M.~Virtanen
  \thanks{
    Department of Mathematics, University of Leicester, Leicester LE1 7RH, England UK.
    \email{jv77@le.ac.uk}
  }
}
\date{\today}                                           % Activate to display a given date or no date
\begin{document}
\maketitle

\begin{abstract}
We consider second order explicit and implicit two-step time-discrete schemes for wave-type equations. We derive optimal order a posteriori estimates 
controlling the time discretization error.  
Our analysis, has been motivated by the need to provide a posteriori estimates for the popular \emph{leap-frog} method (also known as \emph{Verlet}'s method in molecular dynamics literature); it is  extended, however, to general cosine-type second order methods.
The estimators are based on a novel reconstruction of the time-dependent component of the approximation. Numerical experiments confirm similarity of convergence rates of the proposed estimators and of the theoretical convergence rate of the true error. 

\end{abstract}

\section{Introduction}

This work is concerned with second order explicit and implicit two-step time-discrete schemes for wave-type equations.  Our objective is to derive optimal order a posteriori estimates 
controlling the time-discretization error. To the best of our knowledge, error control for wave equations, discretized by popular methods is limited so far to 
first order schemes \cite {bernardi_suli, GLM2013}.  
Despite the importance of such wave-type problems, the lack of error control for  time-discretizations used extensively in applications is probably due to 
the two-step character of these methods, and the associated technical issues. 
Our analysis, has been motivated by the need to provide a posteriori estimates for the \emph{leap-frog} method or, as often termed, \emph{Verlet}'s method in the molecular dynamics literature. It extends, however to general 
cosine-type second order methods \cite{baker-dougalis-serbin,baker-dougalis-serbin2}. \\

%A posteriori error bounds are well developed for stationary boundary
%value problems \citep[e.g.,][and
%the references therein]{verfurth, ainsworth_oden,
%  babuska_strouboulis, cc-zz, dorfler, stevenson, cascon_nochetto}.
%  
  Adaptivity and  a posteriori error control  for
parabolic problems have been developed in, e.g.,  \cite{eriksson-johnson:2,
  verfurth_parabolic,picasso,houston-suli_lagrange,makr-noch,bergam-bernardi-mghazli:05,
  bernardi-verfuerth:05,lakkis-makr}. In particular, as far as time discretization is concerned, 
 all implicit one-step methods can be treated within the framework developed in 
 \cite{amn, makr_noch_diss, AMN2009,AMN2011, LM2013}.
Although some of these results apply (directly or after appropriate modifications) to the wave equation also, when written as a first order system and discretised by implicit Runge-Kutta or Galerkin schemes, 
this framework does not cover popular two-step implicit or explicit time-discretisation methods.
 The recent results  in  \cite{bernardi_suli, GLM2013} cover only first order time discrete schemes; see also \cite{adjerid}
for certain estimators to standard
implicit time-stepping finite element approximations of the wave equation.
For earlier works on adaptivity for wave equations from various perspectives we refer, e.g., to 
 \cite{johnson_wave, bangerth-rannacher:01,suli,suli2}.

%\section{Preliminaries}\label{prelim}
%\noindent
\subsubsection*{Model problem and notation.}\label{mod}
Let $(H,\ltwoin{\cdot}{\cdot})$ be a Hilbert space and $\mA:[0,T]\to D(\mA)$, positive definite, self-adjoint, linear operator on $D(\mA)$, the domain of $\mA$, which is assumed to be dense in $H$. For time $t\in(0,T]$, we consider the linear second order hyperbolic problem: find $u:[0,T]\to D(\mA)$, such that
\begin{equation}\label{pde}
\begin{aligned}
 u''(t)+\mA u(t) &=  f(t)\quad\text{for } 0<t \le T,\\
u(0) &= u_0,\\
u'(0) &= v_0,
 \end{aligned}
\end{equation}
where $f:[0,T]\to H$, $u_0,v_0\in H$.  
%
%
%****** methods for it appears that 
%The class of second degree \emph{cosine methods}, also known as the \emph{Newmark family} of numerical methods for wave problems is used extensively in practical computations. This family includes the very popular \emph{leap-frog method} for a particular choice of method parameters.
%
%This work is concerned with the derivation of optimal order a posteriori error bounds for time-discretizations based on second degree cosine methods. A posteriori bounds for the wave equation have been considered in ...
% 
%{\tt We consider the time-discrete leap-frog method for the wave equation for the moment, but we could do all Newmark methods.}

\subsubsection*{Leap-frog time-discrete schemes.} We shall be concerned with the popular leap-frog time-discrete scheme for
(\ref{pde}). We consider a subdivision of the time interval $(0,T]$ into disjoint subintervals $(t^{n},t^{n+1}]$, $n=0,\dots, N-1$, with $t^0=0$ and $t^N=T$, and we define $k _n:=t^{n+1}-t^{n}$, the time-step. For simplicity of the presentation, we shall assume that $k _n =k$ is constant, although this is not a restriction of the analysis below. Despite being two-step, 
the schemes considered herein can be formulated for variable time steps also, with their consistency and stability properties then being influenced accordingly, cf. \cite{Skeel93, CalvoSanzSerna93}; the study of such extensions is out of the scope of this work. We shall use the notation $t^{n+1/2}:=(t^{n+1}+t^{n})/2$. 
%
%
%{\tt Assume for the moment that $k _n=k $, i.e., a uniform mesh. For non-uniform mesh see later section.} 
%
%\emph{Time discretization methods.}

The time-discrete leap-frog scheme (or Verlet's method in the terminology of initial value problems or of molecular dynamics) for the wave problem (\ref{pde}) is defined by finding approximations 
$U^{n+1}\in D(\mA)$ of the exact values $u^{n+1}:=u(t^{n+1})$, such that:
\begin{equation}\label{lf}
\partial^2 U^{n+1}+\mA U^{n}  =   f^{n},\quad n=1,\dots, N-1,
\end{equation}
where $f^{n}  :=  f(t^{n})\in H$,
\begin{equation}
\partial^2 U^{n+1}   :=  \frac{\partial U^{n+1}-\partial U^{n}}{k  }  =  \frac{ U^{n+1}- 2 U^{n} + U^{n-1}  } {k ^2}
,
\end{equation}
with
\begin{equation*}
\partial U^{n+1}  :=  \frac{U^{n+1}- U^{n}}{k },
\end{equation*}
assuming knowledge of $U^0$ and $U^1$. We set $U^0:= u_0$ and we define $U^1$ by
\begin{equation}\label{fem_fd_neq1}
 \frac{\partial U^{1}-v_0}{k }+\ha \mA U^{0}  =  \ha f^{0},
\end{equation}
where $f^{0}  :=  f(0)$ and $\partial U^1  :=  (U^1-U^0)/k $.
This is a widely used and remarkable method in many ways: it 
is the only two-step explicit scheme for second order problems which is second order accurate, it has 
important conservation and  geometric properties, as it is symplectic, and it is very natural and simple to formulate and implement. 
We refer to the review article \cite{HLW_Verlet2003} for a thorough discussion. 
Explicit schemes, such as \eqref{lf} are suited for the discretization of wave-type partial differential equations, since 
their implementation requires a mild CFL-type condition of the form $k / h \leq C$  (in contrast to parabolic problems,) where $h$ stands for the
space-discretization parameter.% and, thus, they are very popular among applied scientists.  

%\noindent
\subsubsection*{Cosine methods.}  The leap-frog scheme is a member of a general class of two-step methods for second order evolution problems, which are based on the approximation of cosine and are used extensively in practical computations. 
In a two-step cosine method, for  $n=1,\dots, N-1,$ we seek approximations $U^{n+1}$ such that
\begin{equation}\label{cos}
\partial^2 U^{n+1}+\big [ \, q_1\, \mA U^{n+1} -2 p_1  \mA U^{n} + q_1\, \mA U^{n-1}\, \big ]=   \big [ \, q_1\, f^{n+1} -2 p_1  f ^{n} + q_1\, f ^{n-1}\, \big ],\
\end{equation}
where we assume that $p_1= q_1 - \frac 1 2 $ for second order accuracy; we refer to \cite{baker-dougalis-serbin,baker-dougalis-serbin2} for a detailed discussion and analysis of general multi-step cosine schemes. In this case, the rational function $r(x) = (1 + p_1 x^2)/(1 + q_1 x^2)$ is a second order approximation to the
cosine, in the sense that for $|x|$ sufficiently small, 
\begin{equation}\label{cos_2}
| r(x) - \cos (x) | \leq C \, x^4\, . 
\end{equation}
When $q_1 =0$ the above methods are explicit, and the condition $p_1= q_1 - \frac 1 2 $ implies that the only  explicit second order
member of this family is the leap-frog  method \eqref{lf}.

In this work, we derive a posteriori error bounds in the
$\leb\infty$-in-time-energy-in-space norm of the error. 
 The derived bounds are of optimal order, i.e., of the
same order as the error (which is known to be second order) for the class of the schemes considered \cite{baker-dougalis-serbin,baker-dougalis-serbin2}.
This is verified by the numerical experiments presented herein. 
 Our approach is based on the following ingredients:
first, we rewrite the schemes as one-step system on staggered time grids. In turn, this 
can be seen as a second order perturbation of the staggered midpoint method. Further, 
  by introducing appropriate interpolants, we arrive to a form which 
can be viewed as perturbation of \eqref{pde} written as a first-order system. Finally, we employ an adaptation of the time reconstruction from \cite{amn},
 yielding the desired a posteriori error estimates. An interesting observation is that our estimates hold \emph{without} any additional
 time-step assumption, which at the fully discrete level would correspond to a CFL-type restriction. Thus, in a posteriori analysis, 
 standard stability considerations of time discretisation schemes might influence the behaviour of the estimator, but are not explicitly required; the 
 possible instability is sufficiently reflected by the behaviour of the estimator, see Section 4.
Although not done here, by employing space reconstruction techniques it would be possible to 
   derive error estimates for fully discrete schemes in various norms, using ideas from \cite{makr-noch,lakkis-makr, GLM2013}.  
 
The remaining of this work is organized as follows. In
\S\ref{reformulation} we reformulate the numerical methods
appropriately -- this is a crucial step in our approach. We start with
the leap-frog method and we continue with providing two alternative
reformulations of general cosine methods. In \S\ref{apost_bound} we
introduce appropriate time reconstructions and we derive the error
bounds.  In \S\ref{sec:numerics} we present detailed numerical
experiments which yield experimental orders of convergence for the
estimators that are the same with those of the actual error. Finally,
in \S\ref{finalremark} we draw some concluding remarks.
%%%%%%%%%%%%%%%%%%%%%%%%%%%%%%%%%%%%%%%%%%%%%%%%%%%%%%%%%%%%%%%%%%%%%%%%
\clearpage
\section{Reformulation of the methods}\label{reformulation}
 It will be useful for the analysis to reformulate the methods as a system in two staggered grids. 

%We shall derive an a posteriori error bound for the quantity $ \linf{u-U}{0,T;\leb2\Omega}$.
\subsection{Leap-frog\/}
Starting with the leap-frog method, 
we introduce the auxiliary variable
\begin{equation}
V^{n+1/2}   :=    \partial U^{n+1},
\end{equation}
for $n=0,1,\dots,N-1$, 
and we set $V^{-1/2}:=2v_0-V^{1/2}$. (Note that, then, 
$v_0=(V^{-1/2}+V^{1/2})/2$.) 
Also, we define $U^{-1}:=U^0-k V^{-1/2}$ and we observe that we have 
\[
v_0= \frac{U^1-U^{-1}}{2k}.
\]

Further, we introduce the notation
\begin{equation}
\partial V^{n+1/2}:=\frac{V^{n+1/2}-V^{n-1/2}}{k },\quad n=0,1,\dots,N-1,
\end{equation}
noting that the identity $\partial V^{1/2} = 2(\partial U^1 - v_0)/k  $ also holds.

We can now write the method (\ref{lf}) as a system in the staggered form considered in  \cite{HLW_Verlet2003}:
\begin{equation}\label{lf_pw_two}
\begin{split}
&\partial U^{n+1}-V^{n+1/2}  =  0,\\
\quad 
&\partial V^{n+1/2}+\mA U^{n}  =  f^{n} \, , 
\end{split}
\end{equation}
for $n=0,1,\dots,N-1$. 

Next, our goal is to recast  \eqref{lf_pw_two} using globally defined piecewise linear functions. 
We define $U:[-k  ,T]\to D(\mA)$ to be the piecewise linear interpolant of the sequence $\{U^n\}_{n=-1}^N$, at the points $\{t^n\}_{n=-1}^N$, with $t^{-1}:=-k$. %and we set $U^{n+1/2}:=U(t^{n+1/2})$, $n=0,\dots,N-1$, with $t_{-1/2}:=-t^{-1/2}$. 
In addition, let $V:[-k  /2,t^{N-1/2}]\to D(\mA)$  be the piecewise linear interpolant of the sequence $\{V^{n+1/2}\}_{n=-1}^{N-1}$, at the points $\{t^{n+1/2}\}_{n=-1}^{N-1}.$ 
%and we set $V^{n}:=V(t^{n})$, $n=0,\dots,N-1$. Then, by construction we have the (trivial) identities
Using the notation
\begin{equation}\label{lf_not_2}
\begin{split}
U^{n+1/2}:=&\ U(t^{n+1/2}),\\ %\qquad \text{where} t_{-1/2}:=-t^{-1/2},\\
\quad 
V^{n}:=&\ V(t^{n}) , \qquad n=0,\dots,N-1\, .
\end{split}
\end{equation}
we, then, have 
\begin{equation}\label{identities}
U^{n+1/2} =  \ha (U^{n+1}+U^{n}),\quad V^{n}    =  \ha (V^{n+1/2}+V^{n-1/2}) ,
\end{equation}
for  $n=0,1,\dots,N-1$. 

Hence, in view of (\ref{identities}), (\ref{lf_pw_two}) implies
\begin{equation}\label{lf_pw_three}
\begin{aligned}
\partial U^{n+1}- \ha (V^{n+1}+V^{n})  &=-\frac{1}{4} (V^{n+3/2}-2V^{n+1/2}+V^{n-1/2}),
 \\
\partial V^{n+1/2}+\ha\mA (U^{n+1/2}+U^{n-1/2}) &= f^{n}+\frac{1}{4}\mA (U^{n+1}-2 U^{n}+U^{n-1}),
\end{aligned}
\end{equation}
for $n=0,\dots,N-1$.
Upon defining the piecewise constant residuals 
\[
R_U(t)|_{(t^{n-1/2},t^{n+1/2}]}
\equiv R_U^n 
:= \frac{1}{4}\mA (U^{n+1}-2 U^{n}+U^{n-1}),
\]
\[
  R_V(t)|_{(t^{n},t^{n+1}]}
  \equiv R_V^{n+1/2}
  := -\frac{1}{4} (V^{n+3/2}-2V^{n+1/2}+V^{n-1/2})\, ,
  \]
it is easy to check that, given that the leap-frog method is second order (in both $U^n$ and $V^{n+1/2}$), we have $R_U^n =O(k  ^2)$ and $R_V^{n+1/2}=O(k  ^2).$ Hence, \eqref{lf_pw_three}
can be viewed as a second order perturbation of the staggered mid-point method for \eqref{pde} written as first order system 
\begin{equation}\label{wave_syst}
\begin{aligned}
&u'- v  =0\, ,\\
&v'+\mA u =  f\, .
\end{aligned}
\end{equation}
In what follows, it will be useful to rewrite \eqref{lf_pw_three}
as a perturbation of the continuous system \eqref{wave_syst}. To this end, we introduce two time interpolants onto  piecewise linear functions defined on the staggered 
grids:  
we define $U_1:[0,T]\to D(\mA)$ to be the piecewise linear interpolant of the sequence $\{U^{n+1/2}\}_{n=-1}^{N-1}$ and $V_1:[0,t^{N-1}]\to D(\mA)$ to be the piecewise linear interpolant of the sequence $\{V^{n}\}_{n=0}^{N-1}$. 
Then, (\ref{lf_pw_three}) can be written as
\begin{equation}\label{lf_pw_four}
\begin{aligned}
U'- I_0 V_1  &=  R_V,\\
V'+\mA \tilde{I}_0 U_1&= \tilde{I}_0 f  +R_U,
\end{aligned}
\end{equation}
where we define the interpolators
\begin{equation}
  \label{eqn:def:time-interpolators}
  \begin{split}
    \tilde{I}_0:
    &
    \text{ piecewise constant midpoint interpolator on }
    \{(t^{n-1/2},t^{n+1/2}]\}_{n=0}^{N-1}, 
    \\
    I_0
    &
    \text{ piecewise constant midpoint interpolator on }\{(t^{n-1},t^{n}]\}_{n=1}^{N-1}.
  \end{split}
\end{equation}
This formulation will be the starting point of our analysis in the next section. 
%%%%%%%%%%%%%%%%%%%%%%%%%%%%%%%%%%%%%%%%%%%%%%%%%%%%%%%%%%%%%%%%%%%%%%%%
\subsection{Cosine methods: Formulation 1} We shall see that cosine methods \eqref{cos} can be  reformulated  in a similar way as a staggered system. 
As in the leap-frog case we introduce the auxiliary variable
\begin{equation}
V^{n+1/2}   :=    \partial U^{n+1},
\end{equation}
and we let
\begin{equation}
\partial V^{n+1/2}:=\frac{V^{n+1/2}-V^{n-1/2}}{k },\quad n=0,1,\dots,N-1,
\end{equation}
Then the methods  (\ref{cos})  can be rewritten in  system form:
\begin{equation}\label{cos_pw_two}
\begin{split}
&\partial U^{n+1}-V^{n+1/2}  =\  0,\\
\quad 
&\partial V^{n+1/2}+\big [ \, q_1\, \mA U^{n+1} -2 p_1  \mA U^{n} + q_1\, \mA U^{n-1}\, \big ] \\
 = &  \big [ \, q_1\, f^{n+1} -2 p_1  f ^{n} + q_1\, f ^{n-1}\, \big ]\, , 
\end{split}
\end{equation}
for $n=0,1,\dots,N-1$. 
%
%\begin{equation}\label{cos}
%\partial^2 U^{n+1}+\big [ \, q_1\, \mA U^{n+1} -2 p_1  \mA U^{n} + q_1\, \mA U^{n-1}\, \big ]=   \big [ \, q_1\, f^{n+1} -2 p_1  f ^{n} + q_1\, f ^{n-1}\, \big ],\
%\end{equation}
%where we assume that $p_1= q_1 - \frac 1 2 $ for second order accuracy. In this case the rational function $r(x) = \frac{ 1 + p_1 x^2 } {1 + q_1 x^2}$
%
Using the same notation and conventions as in the leap-frog case, 
we observe, respectively, 
\begin{equation}\label{cos_pw_3}
\begin{split}
\big [ \, q_1\,  &\mA U^{n+1} - 2 p_1  \mA U^{n} + q_1\, \mA U^{n-1}\, \big ]  = \ha\mA (U^{n+1/2}+U^{n-1/2})  \\
& - \ha\mA (U^{n+1/2}+U^{n-1/2})  + \big [ \, q_1\, \mA U^{n+1} -2 p_1  \mA U^{n} + q_1\, \mA U^{n-1}\, \big ] \\
   = & \ha\mA (U^{n+1/2}+U^{n-1/2})   \\
& - \ha\mA (U^{n+1/2}+U^{n-1/2})  + \mA U^n + \big [ \, q_1\, \mA U^{n+1} -2 q_1  \mA U^{n} + q_1\, \mA U^{n-1}\, \big ] \\
 = & \ha\mA (U^{n+1/2}+U^{n-1/2})   \\
& -\frac 1 4  \big [   \mA U^{n+1} - 2    \mA U^{n}  +   \, \mA U^{n-1}\, \big ] + \big [ \, q_1\, \mA U^{n+1} -2 q_1  \mA U^{n} + q_1\, \mA U^{n-1}\, \big ] \\
 = & \ha\mA (U^{n+1/2}+U^{n-1/2})    
  -\frac {(1- 4 q_1)} 4  \big [   \mA U^{n+1} - 2    \mA U^{n}   +  \, \mA U^{n-1}\, \big ]  \, ,
\end{split}
\end{equation}
where we used the fact that $p_1= q_1 - \frac 1 2 .$ 
Therefore, as before we conclude that 
\begin{equation}\label{cos_pw_three}
\begin{aligned}
\partial U^{n+1}- \ha (V^{n+1}+V^{n})  &=-\frac{1}{4} (V^{n+3/2}-2V^{n+1/2}+V^{n-1/2}),
 \\
\partial V^{n+1/2}+\ha\mA (U^{n+1/2}+U^{n-1/2}) &=\tilde  f^{n}+\frac {(1- 4 q_1)} 4\mA (U^{n+1}-2 U^{n}+U^{n-1}),
\end{aligned}
\end{equation}
where $ \tilde  f^{n} =  \big [ \, q_1\, f^{n+1} -2 p_1  f ^{n} + q_1\, f ^{n-1}\, \big ] $ and 
  $n=0,\dots,N-1$.
Let us now define 
\[
\begin{aligned}
R^{\cos}_U(t)|_{(t^{n-1/2},t^{n+1/2}]}\equiv R_U^{\cos, n }:=& \ \frac {(1- 4 q_1)} 4 \mA (U^{n+1}-2 U^{n}+U^{n-1})  \\ &+ q_1 \,  \big [ \,\, f^{n+1} -2    f ^{n} +  \, f ^{n-1}\, \big ] ,
\end{aligned}
\]
\[
  R^{\cos}_V(t)|_{(t^{n},t^{n+1}]}\equiv  R_V^{\cos, n+1/2}:= -\frac{1}{4} (V^{n+3/2}-2V^{n+1/2}+V^{n-1/2})\, .
  \]
As in the leap-frog case, it is easy to check that given that the method is second order, we have $R_U^{\cos,n} =O(k  ^2)$ and $R_V^{\cos, n+1/2}=O(k  ^2).$ Hence, \eqref{lf_pw_three}
can be seen as a second order perturbation of the staggered mid-point method for \eqref{wave_syst}.   
Further, still using the same notation for time interpolants as the leap-frog case, we obtain 
\begin{equation}\label{cos_pw_four}
\begin{aligned}
U'- I_0 V_1  &=  R^{\cos}_V,\\
V'+\mA \tilde{I}_0 U_1&= {\tilde{I}_0 \tilde f } +R^{\cos}_U\, .
\end{aligned}
\end{equation}
where   $\tilde{I}_0 \tilde f  | _{ (t^{n-1/2},t^{n+1/2}] }   = \tilde{f}^n.$ It is interesting to compare \eqref{cos_pw_four} to \eqref{lf_pw_four}. \\

\subsection{Cosine methods: Formulation 2} We briefly discuss an alternative formulation  of cosine methods. This time we let
\begin{equation}
V^{n+1/2}   :=   \big (\, I + k ^2  q_1 \mA\, \big ) \partial U^{n+1},
\end{equation}
and, as before,
\begin{equation}
\partial V^{n+1/2}:=\frac{V^{n+1/2}-V^{n-1/2}}{k  },\quad n=0,1,\dots,N-1.
\end{equation}
Then, using again  $p_1= q_1 - \frac 1 2 ,$ we rewrite  the methods  (\ref{cos}) as
\begin{equation}\label{cos2_pw_two}
\begin{split}
\partial U^{n+1}-V^{n+1/2}  =&  -  k ^2  q_1 \mA\,  \partial U^{n+1},\\
\quad 
\partial V^{n+1/2}+   \mA U^{n}    =&  \big [ \, q_1\, f^{n+1} -2 p_1  f ^{n} + q_1\, f ^{n-1}\, \big ]\, , 
\end{split}
\end{equation}
for $n=0,1,\dots,N-1$. 
%
%\begin{equation}\label{cos}
%\partial^2 U^{n+1}+\big [ \, q_1\, \mA U^{n+1} -2 p_1  \mA U^{n} + q_1\, \mA U^{n-1}\, \big ]=   \big [ \, q_1\, f^{n+1} -2 p_1  f ^{n} + q_1\, f ^{n-1}\, \big ],\
%\end{equation}
%where we assume that $p_1= q_1 - \frac 1 2 $ for second order accuracy. In this case the rational function $r(x) = \frac{ 1 + p_1 x^2 } {1 + q_1 x^2}$
%
Using the same notation and conventions as in the leap-frog case, 
we finally conclude
 \begin{equation}\label{cos2_pw_three}
\begin{aligned}
&\partial U^{n+1}- \ha (V^{n+1}+V^{n})  =-  k ^2  q_1 \mA\,  \partial U^{n+1}-\frac{1}{4} (V^{n+3/2}-2V^{n+1/2}+V^{n-1/2}),
 \\
&\partial V^{n+1/2}+\ha\mA (U^{n+1/2}+U^{n-1/2}) =\tilde  f^{n}+\frac {1  } 4\mA (U^{n+1}-2 U^{n}+U^{n-1}),
\end{aligned}
\end{equation}
where $ \tilde  f^{n} =  \big [ \, q_1\, f^{n+1} -2 p_1  f ^{n} + q_1\, f ^{n-1}\, \big ] $ and 
  $n=0,\dots,N-1$.
Upon defining the perturbations as 
\[
R^{\cos,2}_U(t)|_{(t^{n-1/2},t^{n+1/2}]}\equiv R_U^{\cos,2, n }:= \frac {1} 4 \mA (U^{n+1}-2 U^{n}+U^{n-1})  + q_1 \,  \big [ \,\, f^{n+1} -2    f ^{n} +  \, f ^{n-1}\, \big ] ,
\]
\[
  R^{\cos,2}_V(t)|_{(t^{n},t^{n+1}]}\equiv  R_V^{\cos,2, n+1/2}:= -  k ^2  q_1 \mA\,  \partial U^{n+1} -\frac{1}{4} (V^{n+3/2}-2V^{n+1/2}+V^{n-1/2})\, ,
  \]
  it is easy to check, again, that the method is second order ($R^{\cos,2, n} _U  =O(k ^2)$ and $R_V^{\cos,2, n+1/2}=O(k ^2),$) and, thus,  \eqref{lf_pw_three}
can be interpreted  as a second order perturbation of the staggered mid-point method for \eqref{wave_syst}.   
Still, using the same notation as before, we have
\begin{equation}\label{cos2_pw_four}
\begin{aligned}
U'- I_0 V_1  &=  R^{\cos,2}_V,\\
V'+\mA \tilde{I}_0 U_1&= \tilde{I}_0  \tilde f  +R^{\cos,2}_U\, .
\end{aligned}
\end{equation}

\section{A posteriori error bounds}\label{apost_bound}
We have seen that all above schemes can be written in the form, %Then, (\ref{lf_pw_three}) can be written as
\begin{equation}\label{lf_pw_four2}
\begin{aligned}
V'+\mA \tilde{I}_0 U_1&= \tilde{I}_0 (f+\rho_U),
\\
U'- I_0 V_1  &=I_0 \rho_V,
\end{aligned}
\end{equation}
where $\tilde{I}_0,I_0$ are defined in \eqref{eqn:def:time-interpolators} and $\tilde{I}_0\rho_U$ equal to $R_U, R^{\cos}_U$, or  $R^{\cos,2}_U $, for the leap-frog, the first and second cosine method formulations, respectively; similarly, $I_0 \rho_V$ is equal to $R_V, R^{\cos}_V$, or  $R^{\cos,2}_V$ for each of the respective 3 formulations; cf.,  \eqref{lf_pw_four}, \eqref{cos_pw_four} and  \eqref{cos2_pw_four}. It is possible, in principle, to consider non-constant $\rho_U, \rho_V$ on each time-step; nevertheless the, easiest to implement, constant ones considered here suffice to deliver optimal estimator convergence rates, as will be highlighted in the numerical experiments below.

\subsection{Reconstructions} We continue by defining appropriate time reconstructions, cf., \cite{amn}.  To this end, 
 on each interval $(t^{n-1/2},t^{n+1/2}]$, for $n=0,\dots,N-1$, we define the reconstruction $\hat{V}$ of $V$ by
\[
\hat{V}(t):=V^{n-1/2}+\int_{t_{n-1/2}}^t(-\mA  U_1+ \tilde{I}_1 f+\rho_U)\, \ud t ,
\]
where  $\tilde{I}_1$ is a piecewise linear interpolant on the mesh $\{(t^{n-1/2},t^{n+1/2}]\}_{n=1}^{N-1}$, such that $\tilde{I}_1f(t^n)=f^n$. We observe that $\hat{V}(t^{n-1/2})=V^{n-1/2}$, and
\[
\hat{V}(t^{n+1/2})=V^{n-1/2}+k (-\mA  U_1(t^n)+ \tilde{I}_1 f(t^n)+\rho_U(t^n))= V^{n+1/2},
\]
using the mid-point rule to evaluate the integral and the first equation in (\ref{lf_pw_three}).

Also, on each interval $(t^{n-1},t^{n}]$, for $n=1,\dots,N-1$, we define the reconstruction $\hat{U}$ of $U$ by
\[
\hat{U}(t):=U^{n-1}+\int_{t_{n-1}}^t(V_1+\rho_V)\,\ud t.
\]
%where $I_1$ is a suitable piecewise linear interpolant on the mesh $\{(t^{n-1},t^{n}]\}_{n=2}^{N-1}$, such that $I_1f(t^n)=f^n$.
Again, we observe that $\hat{U}(t^{n-1})=U^{n-1}$ and that
\[
\hat{U}(t^n)=U^{n-1}+k ( V_1(t^{n-1/2})+\rho_V(t^{n-1/2}))= U^{n},
\]
using the mid-point rule. Notice that each of the above reconstructions is similar in spirit to the Crank-Nicolson reconstruction of \cite{amn}; however,
we note that, although $\hat U , \hat V$ are both globally continuous functions, their derivatives jump alternatingly at the nodes of the staggered grid.

\subsection{Error equation and estimators} Setting $\hat{e}_U:= u - \hat{U}$ and $\hat{e}_V:= u' - \hat{V}$, we deduce
\begin{equation}\label{lf_recon}
\begin{aligned}
\hat{e}_V'+\mA\hat{e}_U & = \mR_1 + \mR_f \\
\hat{e}_U'-\hat{e}_V & = \mR_2 ,
\end{aligned}
\end{equation}
with
\begin{equation}\label{residual}
\begin{aligned}
 \mR_1 & := -\mA (\hat{U}- U_1)- \rho_U,\\
 \mR_2 & := \hat{V}-  V_1- \rho_V,\\
 \mR_f & := f- \tilde{I}_1 f.
\end{aligned}
\end{equation}
 
For $\Phi=(\phi_1,\phi_2)$, $\Psi=(\psi_1,\psi_2)$ $\in D(\mA)\times H$,  we define the bilinear form
\[
\enin{\Phi}{\Psi}:=\ltwoin{ \mA^{1/2} \phi_1}{\mA^{1/2} \psi_1} + \ltwoin{\phi_2}{\psi_2}.
\]
It is evident that $\enin{\cdot}{\cdot}$ is an inner product on $[D(\mA^{1/2})\times H]^2$. This is the standard 
energy inner product and the induced norm, denoted by $\ndg{\cdot}, $ i.e., 
$$\ndg{\Phi}=(\norm{\mA^{1/2}\phi_1}{}^2+\norm{\phi_2}{}^2)^{1/2}\, , $$
is the natural energy norm for \eqref{pde}.

Then, the a posteriori error estimates will follow by applying standard energy arguments to the error equation \eqref{lf_recon}. More specifically, in view of (\ref{lf_recon}), we have
\[
\begin{aligned}
\ha \frac{\ud}{\ud t}\ndg{(\hat{e}_U,\hat{e}_V)}^2 
& =
\enin{(\hat{e}_U',\hat{e}_V')}{(\hat{e}_U,\hat{e}_V)} \\
& =
\ltwoin{\mA \hat{e}_U'}{\hat{e}_U}+\ltwoin{\hat{e}_V'}{\hat{e}_V}\\
&=
\ltwoin{\mA \hat{e}_V}{\hat{e}_U}+\ltwoin{\mA \mR_2}{\hat{e}_U}-\ltwoin{\mA\hat{e}_U}{\hat{e}_V}+
\ltwoin{\mR_1}{\hat{e}_V}+\ltwoin{\mR_f}{\hat{e}_V}\\
&=
\ltwoin{\mA \mR_2}{\hat{e}_U}+\ltwoin{\mR_1}{\hat{e}_V}+\ltwoin{\mR_f}{\hat{e}_V},
\end{aligned}
\]
using the self-adjointness of $\mA$. Hence, using the Cauchy-Schwarz inequality, we arrive to
\[
\begin{aligned}
\ha \frac{\ud}{\ud t}\ndg{(\hat{e}_U,\hat{e}_V)}^2 
& \le \ndg{(\mR_2,\mR_1+\mR_f)}\ndg{(\hat{e}_U,\hat{e}_V)}.
\end{aligned}
\]
Integrating between $0$ and $\tau$, with $0\le \tau\le t^N$ such that 
\[
\ndg{(\hat{e}_U,\hat{e}_V)(\tau)}{}=\sup_{t\in[0, t^N]}\ndg{(\hat{e}_U,\hat{e}_V)(t)},
\]
we arrive to
\[
\ha \ndg{(\hat{e}_U,\hat{e}_V)(\tau)}^2 
 \le  \ha \ndg{(\hat{e}_U,\hat{e}_V)(0)}^2
+ \ndg{(\hat{e}_U,\hat{e}_V)(\tau)}\int_0^{\tau}\ndg{(\mR_2,\mR_1+\mR_f)}\ud t,
\]
which implies
\[
 \ndg{(\hat{e}_U,\hat{e}_V)(\tau)}^2 
 \le  2 \ndg{(\hat{e}_U,\hat{e}_V)(0)}^2
+ 4\Big(\int_0^{\tau}\ndg{(\mR_2,\mR_1+\mR_f)}\ud t\Big)^2.
\]
This already gives the following a posteriori bound.
\begin{theorem}  \label{MT}
Let $u$ be the solution of \eqref{pde},  $\hat{e}_U:= u - \hat{U}$ and $\hat{e}_V:= u' - \hat{V}.$  Then, the following a posteriori error estimate holds  
\[
\sup_{t\in[0,t^N]} \ndg{(\hat{e}_U,\hat{e}_V)(t)}^2 
 \le  2 \ndg{(\hat{e}_U,\hat{e}_V)(0)}^2
+ 4\Big(\int_0^{t^N}\ndg{(\mR_2,\mR_1+\mR_f)}\ud t\Big)^2,
\]
where $\mR_2,\mR_1$, and $\mR_f$ are defined in \eqref{residual}. 
\end{theorem}
An immediate Corollary from Theorem \ref{MT} is an posteriori bound for the
error $\sup_{[0,t^N]} \ndg{(u-U, u'-V)}$ which can be trivially deduced
through a triangle inequality.
\begin{remark}
\smallskip
Notice that due to the two-mesh stagerring, the computation of the `last' $V$ used in the above estimate, $V^{N-1/2},$ requires the computation of $U^{N+1}.$
This can be obtained by advancing one more time step in the computation before estimating. 
\end{remark}
%%%%%%%%%%%%%%%%%Numerical Examples
%%%%%%%%%%%%%%%%%%%%%%%%%%%%%%%%%%%%%%%%%%%%%%%%%%%%%%%%%%%%%%%%%%%%%%%%
\section{Numerical Experiments}
\label{sec:numerics}
\subsection{Fully discrete formulation.}
Although the focus of the present work is in time-discretization, we shall introduce a fully discrete version of the time-stepping schemes for the numerical experiments bellow.  To this end, we consider  
%This work is concerned with the derivation and numerical study of a posteriori error estimates in $\lebspace\infty(|||\cdot|||)$-norm ($\lebspace\infty$-energy-norm) for a fully discrete approximation of the problem: 
Let $\dom \subset \mathbb{R}^d$, $d=2,3$ be a domain with boundary $\partial\dom$. We consider the initial-boundary value problem for the wave equation: find $\pdesol\in \lebspace{\infty}(0,T;\sobu10\Omega)$ such that,
\begin{eqnarray} 
\pdesol_{tt} +  \mathcal {A} \pdesol &=& \pdesource \quad \;  \text{ in } \; \dom \times (0,T], \label{modelpde1}\\
\pdesol &= &\initpdesol \quad \text{ in } \; \dom \times \{0\}  \quad  \label{modelpde2}\\
\pdesol_t &= &\initpdevel \quad \text{ in } \; \dom \times \{0\} \quad \text{} \label{modelpde3}\\
\pdesol &= &g \quad \text{ on } \; \partial\dom \times (0,T] ,  
\end{eqnarray}
where, for simplicity, we  take $\mathcal {A}  =-c^2\Delta, c\ne 0,$ and $g\in\sob{1/2}{\partial\Omega}$.  
%
%More specifically, we perfform a numerical study of an a posteriori error estimator for a numerical scheme consisting of continuous Galerkin method in space and leapfrog time-stepping for the problem  (\ref{modelpde1}) - (\ref{modelpde3}). 
Further, for each $n$, we consider the standard, conforming finite element space $S_h^p\subset\sobu10\Omega$, based on a quasiuniform triangulation of $\Omega$ consisting of finite elements of polynomial degree $p$, with $h$ denoting the largest element diameter. Focusing on time-discretization issues, we shall use the same spatial discretization for all time-steps.  The respective discrete spatial operator is denoted by $\mathcal{A}_h^p$.  The fully-discrete leap-frog method is then defined as follows: for each $n=2,\dots,N$, find $ \usoln^{n+1} \in S_h^p$, such that
%#<u_{n+1},v>=<2 u_n - u_{n-1},v> + (\Delta_t^2)(<f(t_n),v>-B(u_n,v))
\begin{equation}\label{eqn:fully-discrete-leapfog}
\begin{aligned}
    \usoln^{n+1} &=   2 \usoln^{n} -   \usoln^{n-1} + \timestep ^2 (  \pdesourcetimeapprox^{n} -\mathcal{A}_h^p\usoln^{n} ) ,
  \end{aligned}
\end{equation}
and $\usoln^{1} \in  S_h^p$ such that
%#<u_1,v>=<u_0,v> + <\Delta_t v_0,v> + (\Delta_t^2 / 2)(B(u_0,v) - <f(0),v>)
\begin{equation}\label{eqn:fully-discrete-leapfog-initial-step}
  \begin{aligned}
 \usoln^{1}&=  \usoln^{0} + \timestep   \velocsoln^0+ \frac{\timestep ^2}{2} ( \mathcal{A}_h^p \usoln^{0} - \pdesourcetimeapprox^{0} ) ;  
  \end{aligned}
\end{equation}
here %$\{ t_n : t_n \in [0,T], \; \forall \; n = 0 \ldots N \}$, $\timestep =t_n-t_{n-1}$ for $n \ge 1 $ and $\timestep_0:=\timestep $, 
 $\pdesourcetimeapprox^n(\cdot):=\Pi \pdesource(\cdot,t_n)$ for each $t_n$, where $\Pi: \lebspace2(\dom) \to S_h^p$ denotes a suitable interpolation/projection operator onto the finite element space $S_h^p$ %, and {\color{red}elsewhere a piecewise $\lebspace2$-projection of polynomial degree $k$ in time MEANING??} 
 of the source function $\pdesource$. We also set $\usoln^0:=\Pi \initpdesol$ and $\velocsoln^0:=\Pi \initpdevel$. Note that $V$ can be calculated as above through $U$, or can be computed as follows:
find  $\velocsoln^{n+1/2} \in S_h^p$ such that 
\begin{equation}\label{eqn:fully-discrete-leapfog-velocity-formulation}
  \begin{aligned}
     \velocsoln^{n+1/2} &=  \velocsoln^{n-1/2} + \timestep (  \pdesourcetimeapprox^{n}-\mathcal{A}_h^p\usoln^{n} )  ; &
  \end{aligned}
\end{equation}
\eqref{eqn:fully-discrete-leapfog-velocity-formulation} can be used to overcome the difficulty of evaluating estimators defined on staggered time mesh and depending on the term $\velocsoln^{n+3/2}$. 

To assess the time-error estimator, we replace $\mathcal{A}$ by its approximation $\mathcal{A}_h^p$ in the a posteriori estimators discussed above, and in the $((\cdot,\cdot))$-inner product and $\enorm{\cdot}$-norm. For brevity, we introduce the notation:
\[
e_R := (\hat{e}_U,\hat{e}_V)\quad\text{and}\quad  e_L := e_R+(U-\hat{U},V-\hat{V}) = (u - U, u_t -V).
\]
The objective is to study the performance of the a posteriori estimator 
 \begin{equation}\label{def:estimator1}
    %\eta_1 := {\left( 2\enorm{e_R(0)}^2 + 4\left( {\left(\int_0^{T} {B(R_2,R_2)}^{1/2} \right)}^2 + {\left(\int_0^{T} \norm{R_1+R_f}{} \right)}^2 \right) \right)}^{1/2}
    \eta_1 := {\left( 2\enorm{e_R(0)}^2 + 4 {\left(\int_0^{T} \enorm{(R_2,R_1+R_f)} \right)}^2 \right)}^{1/2},
\end{equation}
from Theorem \ref{MT}.
%%R1, R2
%\begin{equation}\label{def:R1}
%    R_1 := -\mathcal{A}\left(\urecsoln - \lineintopIhalfstep \usolnhalfstep \right) - \rho_U, 
%\end{equation}
%and
%\begin{equation}\label{def:R2}
%    R_2 := \recvelocsoln - \lineintopIwholestep \velocsolnwholestep - \rho_V
%\end{equation}
%and
%\begin{equation}\label{def:Rf}
%    R_f := \pdesource - \lineintopIhalfstep \pdesource\, .
%\end{equation}

%\begin{remark}
%Note that 
%\begin{eqnarray}
%{B(R_2,R_2)}^{1/2} &\le & \enorm{(R_2,R_1+R_f)} \; , \\
%\norm{R_1+R_f}{} & \le & \enorm{(R_2,R_1+R_f)}  \; \text{and} \\
%\enorm{(R_2,R_1+R_f)} & \le & {B(R_2,R_2)}^{1/2}+\norm{R_1+R_f}{} 
%\end{eqnarray}
%which is to say that the bound in \ref{theorem:apost_hyperbolic} is of the same order as in \cite{glm_wave:2011}.
%\end{remark}

\subsection{Specific tests}\label{numerics}

For $t\in [0,1]$ and $\dom := (0,1)^2$, we consider the model problem (\ref{modelpde1}) - (\ref{modelpde3}) in the following set up:
$  \mathcal{A} := - c^2 \Delta$ and $ \pdesource=0$,
for $c > 0$ constant. Then, the exact solution $u$ of problem (\ref{modelpde1}) - (\ref{modelpde3}) is given by:
\begin{equation}\label{def:exactsoln}
    u(t,x) := \sum_{j,k=1}^{\infty} \sin(k \pi x)\sin(j \pi x) \left( \alpha_{k,j}\cos(\xi_{k,j} \pi t) +\beta_{k,j}\sin(\xi_{k,j} \pi t)\right)
\end{equation}
where $\alpha_{k,j}>0$, $\beta_{k,j}>0$ and $\xi_{k,j}:=c \sqrt{k^2+j^2}$.
To illustrate the estimator's behaviour, we have chosen the following sets of parameters as numerical examples:
\begin{eqnarray}
 \text{ } & & \left\{
  \begin{aligned}
  \label{eqn:numerics:example1} 
  \text{ } &c:=1.0,\\
  \text{ } &\alpha_{1,1}=\beta_{1,1}=15.0, \beta_{k,j}=\alpha_{k,j}=0 \; \text{for }\;k,j \neq 1, \\
%  \text{ } &\text{ sequence of mesh sizes:} \{1/2 \text{(cyan)},1/4 \text{(green)},1/(4\sqrt{2}) \text{(yellow)},1/8 \text{(red)},1/10 \text{(purple)}\}
  \end{aligned}
  \right. \\
  \text{ } & & \left\{
  \begin{aligned}
  \label{eqn:numerics:example2}
  \text{ } &c:=1.0,\\
  \text{ } &\alpha_{3,3}=\beta_{3,3}=1.0, \beta_{k,j}=\alpha_{k,j}=0 \; \text{for }\;k,j \neq 3, \\
  %\text{ } &\text{ sequence of mesh sizes:} (1/2,1/4,1/(4\sqrt{2}),1/8,1/10) \;\text{for}\; p=3 \\
 % \text{ } &\text{ sequence of mesh sizes:} \{1/(4\sqrt{2}) \text{(cyan)},1/8 \text{(green)},1/10 \text{(yellow)},1/12 \text{(red)},1/14 \text{(purple)}\}
  \end{aligned}
  \right. \\
\text{ } & &   \left\{
  \begin{aligned}
  \label{eqn:numerics:example3}
  \text{ } &c:=5.0,\\
  \text{ } &\alpha_{1,1}=\beta_{1,1}=15.0, \beta_{k,j}=\alpha_{k,j}=0 \; \text{for }\;k,j \neq 1, \\
 % \text{ } &\text{ sequence of mesh sizes:} \{1/2 \text{(cyan)},1/4 \text{(green)},1/(4\sqrt{2}) \text{(yellow)},1/8 \text{(red)},1/10 \text{(purple)}\}
  \end{aligned}
  \right. 
\end{eqnarray}

The solutions of (\ref{eqn:numerics:example1}) - (\ref{eqn:numerics:example3}) are all smooth, but (\ref{eqn:numerics:example3}) oscillates much faster temporally, while (\ref{eqn:numerics:example2}) has greater space-dependence of the error. In the numerical experiments, the C++ library FEniCS/dolfin 1.2.0 and PETSc/SuperLU were used for the finite element formulation and the linear algebra implementation. For each of the examples, we compute the solution of (\ref{eqn:fully-discrete-leapfog}) using finite element spaces of polynomial degree $\polydegree=2$, and time step size $\timestep = C h^r/(p+1)^2$, $r=1,2$, for some constant $C>0$, with $h>0$ denoting the diameter of the largest element in the mesh associated with $S_h^p$. The sequences of meshsizes considered (with respective colouring in the figures below) are: $h=1/2$ (cyan),$1/4$ (green), $1/(4\sqrt{2})$ (yellow), $1/8$ (red), $1/10$ (purple), for Examples (\ref{eqn:numerics:example1}) and (\ref{eqn:numerics:example3}), and  $h=1/(4\sqrt{2})$ (cyan), $1/8$ (green), $1/10$ (yellow), $1/12$ (red), $1/14$ (purple) for Example (\ref{eqn:numerics:example2}). Note that the CFL-condition required by the leap-frog method, is satisfied by the restriction on the time step
\begin{equation}
  \label{eqn:cfl-condition}
    \timestep \le C \frac{h}{(\polydegree+1)^2},
\end{equation}
for sufficiently small $C>0$ constant, which couples the temporal and spatial discretization sizes. We also include an experiment to highlight the behaviour of the estimator when 
the CFL condition is violated.
% The quantities of interest is the error by linearly interpolated sequence of discrete solutions and velocities; ${\{\usoln^n\}}_{n=1,\dots, N}$ and ${\{\velocsoln^{n+1/2}\}}_{n=0,\dots, N - 1}$ are extended into continuous piecewise linear functions of time by,
%\begin{equation} \label{eqn:fully-discrete-solution-extended}
%\usoln := \sfrac{t - t^{n-1}}{\timestep_n} \usoln^n + (1 - \sfrac{t - t^{n-1}}{\timestep_n}) \usoln^{n-1} \quad \velocsoln := \sfrac{t - t^{n-1/2}}{t^{n+1/2}-t^{n-1/2}} \velocsoln^{n+1/2} + (1 - \sfrac{t - t^{n-1/2}}{t^{n+1/2}-t^{n-1/2}}) \velocsoln^{n-1/2}
%\end{equation}
%respectively and 

We monitor the evolution of the values and the experimental order of
convergence of the estimator $\eta_1$, and the errors $e_R$ and $e_L$, as well as the effectivity
index over time on a sequence of uniformly refined meshes with mesh sizes given as per each example. We also monitor the energy of the reconstructed solution:
\begin{equation} \label{eqn:fully-discrete-solution-reconstruction-energy}
E_{\text{reconstruction}}:=(\urecsoln, \recvelocsoln).
\end{equation}

We define \emph{experimental order of convergence} ($\EOC$) of a given sequence of positive quantities $a(i)$ defined on a sequence of meshes of size $h(i)$ by
\begin{equation}
  \EOC( a,i ) = \frac{\log(a(i+1)/a(i))}{\log(h(i+1)/h(i))},
\end{equation}
 the \emph{inverse effectivity index} ,
\begin{equation}
  \IEI(\ltwo{e}{\lebspace\infty(0,t_m;\enorm{\cdot})},\eta_1) = \frac{ \ltwo{e}{\lebspace\infty(0,t_m;\enorm{\cdot})}}{\eta_1} .
\end{equation}
The \IEI has the same information as the (standard) effectivity index and has the advantage of relating 
directly to the inequality appearing in Theorem \ref{MT}. %(\ref{hyperbolic_energy_apost_bound}). 
The results of numerical 
experiments on uniform meshes, depicted in Figures \ref{fig:result1}  and \ref{fig:result2}, indicate that the error estimators 
are reliable and also efficient provided the time steps are kept sufficiently small.  In the last experiment, Figure \ref{CFL-violation},  the behaviour of the estimator is displayed, for the case when the CFL condition \eqref{eqn:cfl-condition} is violated; the estimator remains reliable in this case also: the instability is reflected in the behaviour of estimator, cf., {Figure \ref{CFL-violation}}. Indeed, the inverse effectivity idex $\IEI$ of the estimator oscillates around a constant value for all times in the course of the numerically unstable behaviour.
%%%%%%%%%%%%%%%%%%%%%%%%%%%%%%%%%%%%%%%%%%%%%%%%%%%%%%%%%%%%%%%%%%%%%%%%
\section{Concluding remarks}
\label{finalremark}
An a posteriori error bound for error measured in $\lebspace\infty$ - norm in time and energy norm in space for leap-frog and cosine-type time semi-discretizations for linear second order evolution problems was presented and studied numerically. The estimator was found to be reliable, with the same convergence rate as the theoretical convergence rate of the error. In a fully discrete setting, this estimator corresponds to the control of time discretization error.  %$\BigO{h^\polydegree + \timestep^2}$. 
The estimators were also found to be sharp on uniform meshes provided that the time steps and, thus, the time-dependent part of the error, is kept sufficiently small. Investigation into the suitability of the proposed estimators within an adaptive algorithm remains a future challenge. 

%update refs
\bibliographystyle{siam}

%%%%%%%%%%%%%%%%%%%%%%%%%%%%%%%%%%%%%%%%%%%%%%%%%%%%%%%%%%%%%%%%%%%%%%%%
%% FIGURES START HERE, ALL DUMPED AT END
%%%%%%%%%%%%%%%%%%%%%%%%%%%%%%%%%%%%%%%%%%%%%%%%%%%%%%%%%%%%%%%%%%%%%%%%
\def\vspacingup{-0.5cm}
\def\imagewidth{0.65\textwidth}
\def\imageheight{10.0cm}
%%%%%%%%%%%%%%%%%%%%%%%%%%%%%%%%%%%%%%%%%%%%%%%%%%%%%%%%%%%%%%%%%%%%%%%%
%%\clearpage
%%%%%%%%%%%%%%%%%%%%%%%%%%%%%%%%%%%%%%%%%%%%%%%%%%%%%%%%%%%%%%%%%%%%%%%%
%%%Example 1
%%
%%%%%%%%%%%%%%%%%%%%%%%%%%%%%%%%%%%%%%%%%%%%%%%%%%%%%%%%%%%%%%%%%%%%%%%%
\begin{figure}[htp]
%\begin{figure}
\vspace{\vspacingup}
  \caption{\label{fig:result1} Examples (\ref{eqn:numerics:example1}) and (\ref{eqn:numerics:example2}). }
  \begin{center}
	\vspace{0.0cm}
    \subfigure[{
Example (\ref{eqn:numerics:example1}). Errors, estimator and $\IEI$ are plotted on the top row and $\EOC$s and energy of the reconstructed solution on the bottom row over time ($x$-axis). Results are computed on the sequence of uniform meshes with mesh size $h$, fixed time step $\timestep = 0.4 h / (\polydegree +1)^2 $ and $\polydegree = 2$. The $\IEI$ behaviour indicates that the error is well estimated by the estimator and the convergence rate of the estimator remains near to $\EOC \approx 2$, i.e., to that of the errors $e_L$ and $e_R$. 
    }]{
		\includegraphics[trim=55 15 55 15,clip,width=\imagewidth]{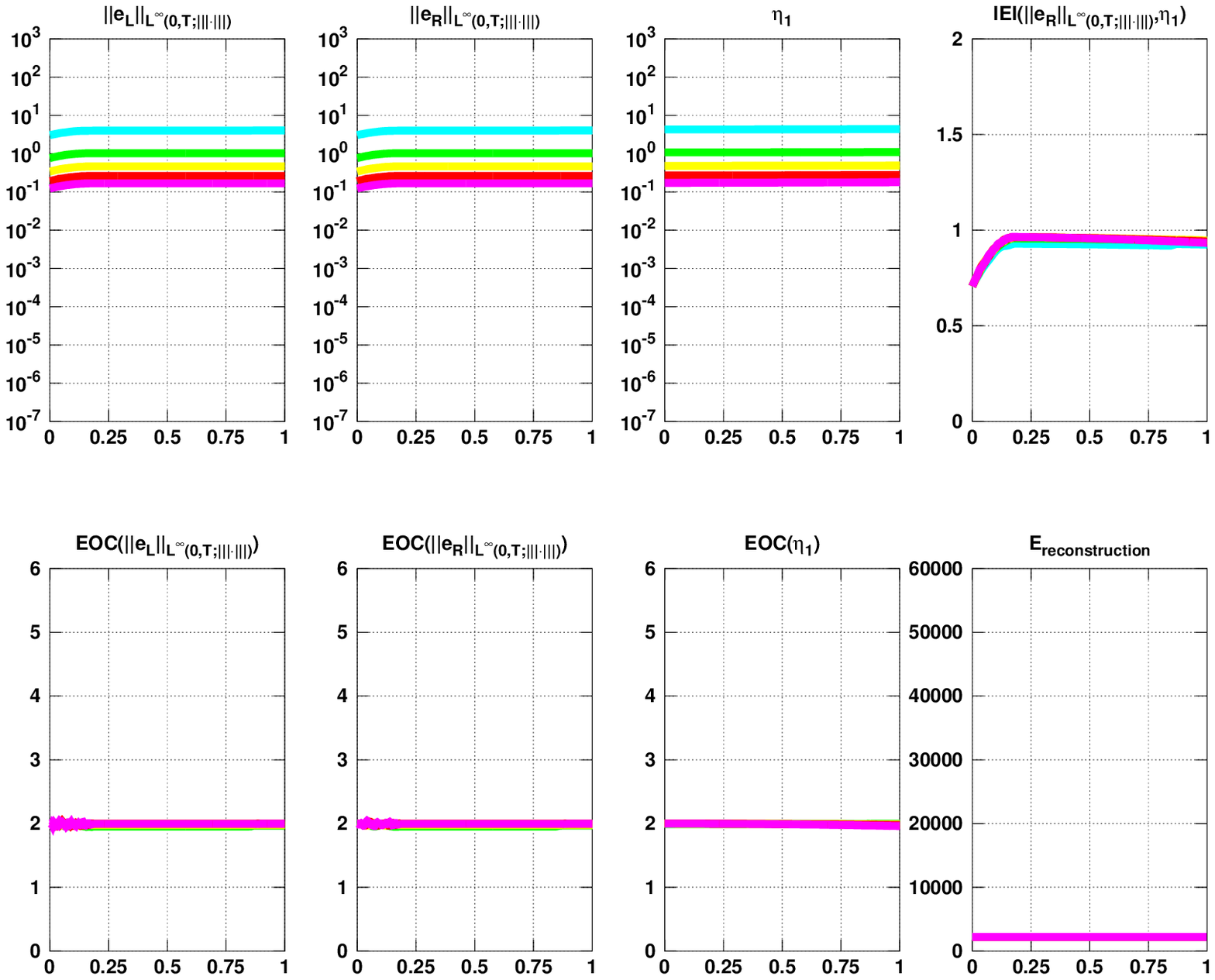}
    }
	\vspace{0.0cm}
    \subfigure[{
Example (\ref{eqn:numerics:example2}). Errors, estimator and $\IEI$ are plotted on the top row and $\EOC$s and energy of the reconstructed solution on the bottom row over time ($x$-axis). Results are computed on the sequence of uniform meshes with mesh size $h$, fixed time step $\timestep = 0.4 h / (\polydegree +1)^2 $ and $\polydegree = 2$. The $\IEI$ behaviour indicates that the error is well estimated by the estimator and the convergence rate of the estimator remains near to $\EOC \approx 2$, i.e., to that of the errors $e_L$ and $e_R$. 
    }]{
		\includegraphics[trim=55 15 55 15,clip,width=\imagewidth]{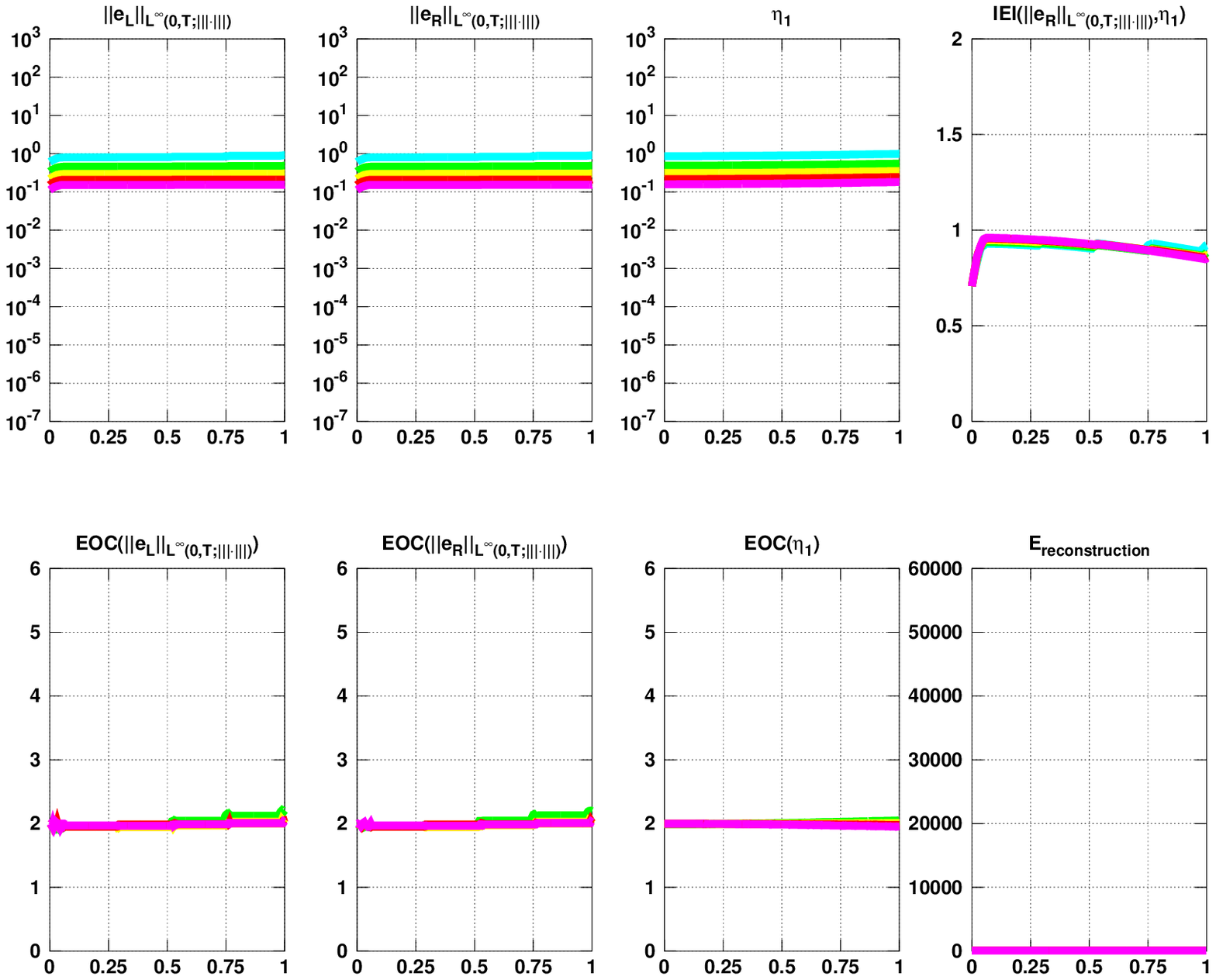}
    }
  \end{center}
\end{figure}
%%%%%%%%%%%%%%%%%%%%%%%%%%%%%%%%%%%%%%%%%%%%%%%%%%%%%%%%%%%%%%%%%%%%%%%%
\begin{figure}[htp]
%\begin{figure}
\vspace{\vspacingup}
  \caption{\label{fig:result2} Example (\ref{eqn:numerics:example3}). }
  \begin{center}
	\vspace{0.0cm}
    \subfigure[{
Example (\ref{eqn:numerics:example3}).  Errors, estimator and $\IEI$ are depicted on the top row and $\EOC$s and energy of the reconstructed solution on the bottom row over time ($x$-axis). Results are computed on the sequence of uniform meshes with mesh size $h$, time-step $\timestep = 0.1 h / (\polydegree +1)^2 $ and $\polydegree = 2$. The $\IEI$ behaviour indicates that the error is overestimated by the estimator. The convergence rate of the estimator is slightly below $2$, probably due to the somewhat coarse time-step for asymptotic convergence, cf., subfigure (b) below.
    }]{
			\includegraphics[trim=55 15 55 15,clip,width=\imagewidth]{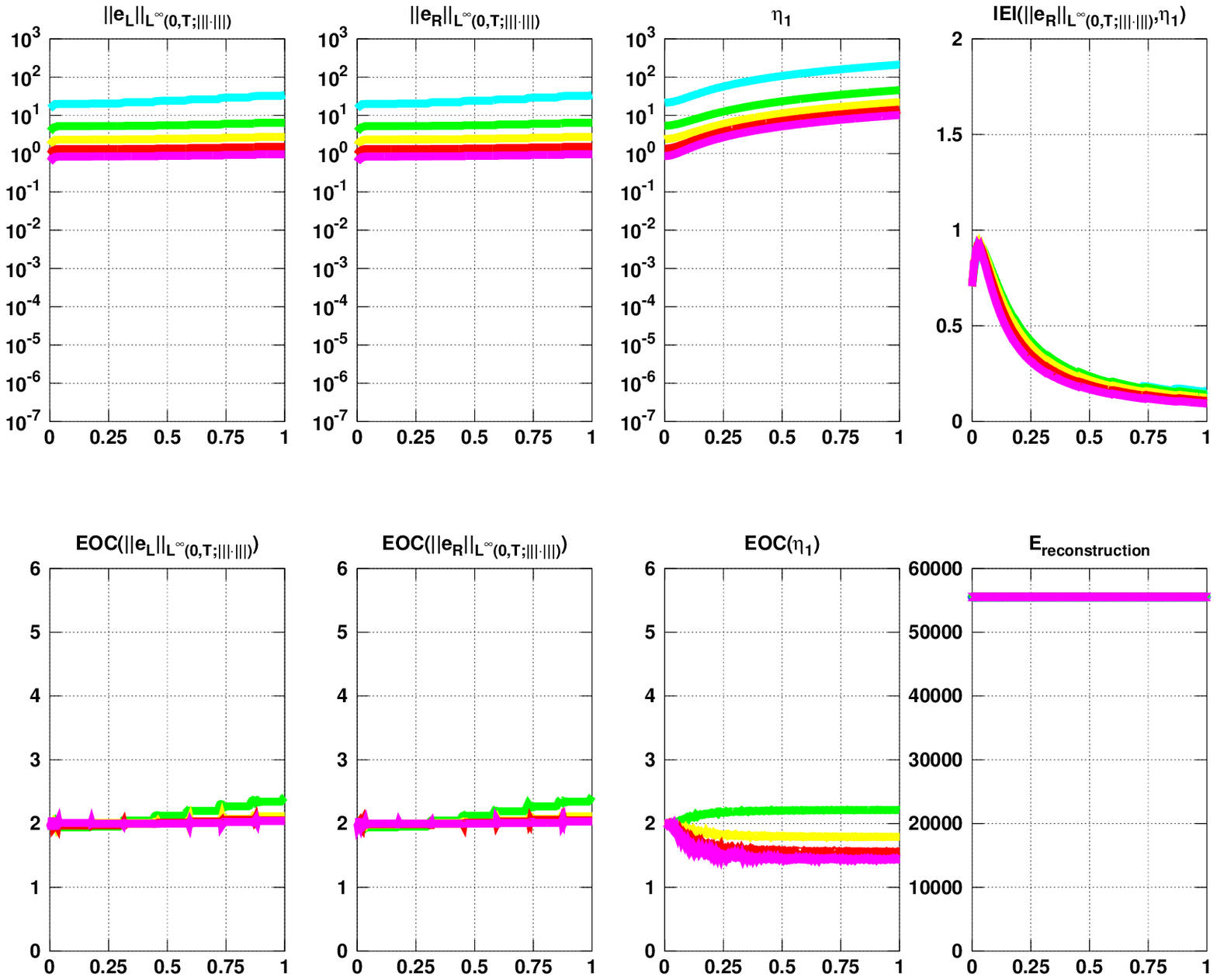}
    }
	\vspace{0.0cm}
    \subfigure[{
Example (\ref{eqn:numerics:example3}). Errors, estimator and $\IEI$ are depicted on the top row and $\EOC$s and energy of the reconstructed solution on the bottom row over time ($x$-axis). Results are computed on the sequence of uniform meshes with mesh size $h$, time-step $\timestep = 0.4 h^2 / (\polydegree +1)^2 $, and $\polydegree = 2$. The effect of using slightly smaller time-step for fine $h$ (resulting from the $h^2$-term) is evident compared to subfigure (a) above, in that the $\EOC \ge 2$ for all times. 
    }]{
			\includegraphics[trim=55 15 55 15,clip,width=\imagewidth]{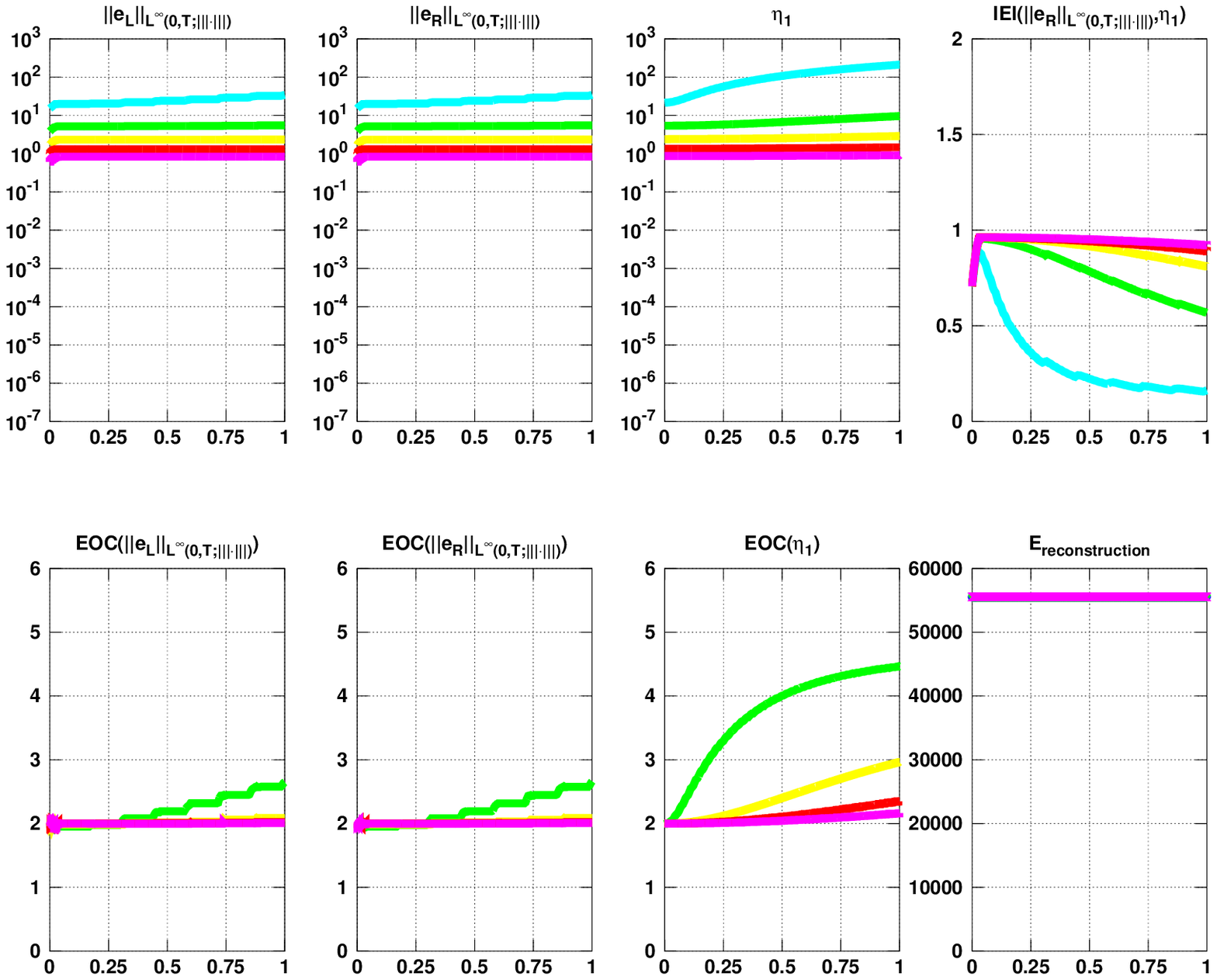}
    }
  \end{center}
\end{figure}

\begin{figure}[]
%\begin{figure}
%\vspace{\vspacingup}
  \caption{\label{fig:result6} Example (\ref{eqn:numerics:example1}), violation of the CFL condition. }
  \begin{center}
	\vspace{-0.1cm}
    \subfigure[{
Example (\ref{eqn:numerics:example1}). Errors, estimator and $\IEI$ are depicted on the top row and $\EOC$s and energy of the reconstructed solution on the bottom row over time ($x$-axis). Results are computed on the sequence of uniform meshes with mesh size $h$ and time step $\timestep = 2.0 h/(\polydegree+1)^2$ and $\polydegree = 3$. The $\IEI$ behaviour indicates that the error is overestimated by the estimator but follows the error behaviour. The method is unstable due to the violation of the CFL condition (cf., Figure \ref{fig:result2}(b) for comparison with a stable approximation.)
    }]{
		\includegraphics[trim=55 15 55 15,clip,width=\imagewidth]{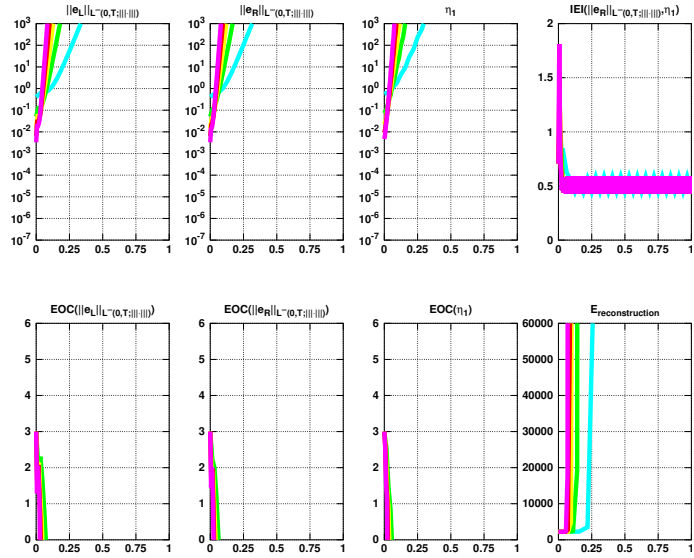}
    }
  \end{center}
  \label{CFL-violation}
\end{figure}
%%%%%%%%%%%%%%%%%%%%%%%%%%%%%%%%%%%%%%%%%%%%%%%%%%%%%%%%%%%%%%%%%%%%%%%%
\end{document}